\documentclass[12pt]{amsart}
\usepackage{amssymb,amsmath}
\numberwithin{equation}{section}

\def\T{\text}
\def\1#1{\overline{#1}}
\def\2#1{\widetilde{#1}}
\def\3#1{\widehat{#1}}
\def\4#1{\mathbb{#1}}
\def\5#1{\frak{#1}}
\def\6#1{{\mathcal{#1}}}
\def\C{{\4C}}
\def\R{{\4R}}

\def\Z{{\4Z}}

\begin{document}
%
\title[Extremal discs and the holomorphic extension]
{Extremal discs and the holomorphic extension from convex hypersurfaces}
\author[ L.~Baracco, A.~Tumanov, G.~Zampieri ]
{Luca Baracco, Alexander Tumanov, Giuseppe Zampieri }
\maketitle
\def\B{\Bbb B}
\def\Giialpha{\mathcal G^{i,i\alpha}}
\def\cn{{\C^n}}
\def\cnn{{\C^{n'}}}
\def\ocn{\2{\C^n}}
\def\ocnn{\2{\C^{n'}}}
\def\const{{\rm const}}
\def\rk{{\rm rank\,}}
\def\id{{\sf id}}
\def\aut{{\sf aut}}
\def\Aut{{\sf Aut}}
\def\CR{{\rm CR}}
\def\GL{{\sf GL}}
\def\Re{{\sf Re}\,}
\def\Im{{\sf Im}\,}
\def\codim{{\rm codim}}
\def\crd{\dim_{{\rm CR}}}
\def\crc{{\rm codim_{CR}}}
\def\phi{\varphi}
\def\eps{\varepsilon}
\def\d{\partial}
\def\a{\alpha}
\def\b{\beta}
\def\g{\gamma}
\def\G{\Gamma}
\def\D{\Delta}
\def\Om{\Omega}
\def\k{\kappa}
\def\l{\lambda}
\def\L{\Lambda}
\def\z{{\bar z}}
\def\w{{\bar w}}
\def\Z{{\1Z}}
\def\t{{\tau}}
\def\th{\theta}
\emergencystretch15pt
\frenchspacing
\newtheorem{Thm}{Theorem}[section]
\newtheorem{Cor}[Thm]{Corollary}
\newtheorem{Pro}[Thm]{Proposition}
\newtheorem{Lem}[Thm]{Lemma}
\theoremstyle{definition}\newtheorem{Def}[Thm]{Definition}
\theoremstyle{remark}
\newtheorem{Rem}[Thm]{Remark}
\newtheorem{Exa}[Thm]{Example}
\newtheorem{Exs}[Thm]{Examples}
\def\Label#1{\label{#1}}
\def\bl{\begin{Lem}}
\def\el{\end{Lem}}
\def\bp{\begin{Pro}}
\def\ep{\end{Pro}}
\def\bt{\begin{Thm}}
\def\et{\end{Thm}}
\def\bc{\begin{Cor}}
\def\ec{\end{Cor}}
\def\bd{\begin{Def}}
\def\ed{\end{Def}}
\def\br{\begin{Rem}}
\def\er{\end{Rem}}
\def\be{\begin{Exa}}
\def\ee{\end{Exa}}
\def\bpf{\begin{proof}}
\def\epf{\end{proof}}
\def\ben{\begin{enumerate}}
\def\een{\end{enumerate}}

\def\simto{\overset\sim\to\to}
\def\1alpha{[\frac1\alpha]}
\def\T{\text}
\def\R{{\Bbb R}}
\def\I{{\Bbb I}}
\def\C{{\Bbb C}}
\def\Z{{\Bbb Z}}
\def\Fialpha{{\mathcal F^{i,\alpha}}}
\def\Fiialpha{{\mathcal F^{i,i\alpha}}}
\def\Figamma{{\mathcal F^{i,\gamma}}}
\def\Real{\Re}

\section{Introduction}
Let $D\subset\C^n$ be a convex domain with smooth boundary $\partial D$ and
let $f$ be a continuous function on $\partial D$.
Suppose for every complex line $L$ the restriction $f|_{L\cap\partial D}$
holomorphically extends into $L\cap D$. Then $f$ extends to $D$ as a
holomorphic function of $n$ variables (Stout \cite{S77}).
The conclusion is still true if instead of the holomorphic extendibility
of $f$ into the sections $L\cap D$, we assume the weaker
Morera condition
\begin{equation}
\Label{morera}
\int_{L\cap\partial D}f\alpha=0
\end{equation}
for every  $(1,0)$-form $\alpha$ with constant coefficients
and every complex line $L$ (Globevnik and Stout \cite{GS}).

The condition of holomorphic extendibility into sections $L\cap D$
and even the Morera condition \eqref{morera} for all lines $L$ seem excessively
strong because it suffices to use only the lines close to the
tangent lines to $\partial D$. Indeed, for simplicity assume $f\in
C^1(\partial D)$. Then the Morera condition for $L$ as $L$
approaches a tangent line $L_0$ at $z_0\in\partial D$ implies that
the $\bar\partial$ derivative of $f$ at $z_0$ along $L_0$ equals
zero. Then $f$ holomorphically extends to $D$ by the classical
Hartogs-Bochner theorem. Therefore of great interest are ``small''
families of lines, for which the result is still true. In
particular, the family of lines should not contain the lines close
to the tangent lines to $\partial D$.

Reducing the family of lines,
Agranovsky and Semenov \cite{AS} show that if $D_2\subset D_1$
are domains in $\C^n$ and $f\in C(\partial D_1)$ holomorphically
extends into sections $L\cap D_1$ by the lines that meet $D_2$,
then $f$ holomorphically extends to $D_1$.
In the case of two concentric balls $D_2\subset D_1$, Rudin \cite {R}
proves that the same conclusion is valid if one only assumes the
extendibility into sections by the lines tangent to $\partial D_2$.
Globevnik \cite{G83} observes that in Rudin's result one only needs
the lines tangent to a sufficiently large open set in $\partial D_2$.

Globevnik and Stout \cite{GS} conjecture that Rudin's result is valid
for every convex domains $D_2\subset\subset D_1$, that is if
$f\in C(\partial D_1)$ holomorphically
extends into sections $L\cap D_1$ by the lines tangent to $\partial D_2$,
then $f$ holomorphically extends to $D_1$.
They also observe in \cite{GS} (see also \cite{BCPZ}) that for $n=2$
in Rudin's result one generally cannot replace the extendibility into
the sections $L\cap D_1$ by the Morera condition \eqref{morera},
that is the latter
suffices unless the ratio $r_1/r_2$ of the radii of the balls
belongs to an exceptional countable set.
For a counterexample in $\C^2$, take $r_1=1$, $r_2=\sqrt{\frac13}$,
and $f=z_1\bar z_2^2$.
However if $n>2$, then Berenstein, Chang, Pascuas, and Zalcman \cite{BCPZ}
show that for the concentric balls the Morera condition for the tangent
lines suffices without exceptions.

Further reduction of the family of lines is possible.
Globevnik \cite{G87} shows that for the unit ball $D\subset\C^2$,
the holomorphic extension property into sections by lines of certain two
parameter family suffices for the holomorphic extendibility into $D$.
The set of lines in his result consists of two disjoint tori.
The second author shows \cite{T97} that for every generating CR manifold
$M\subset\C^n$ of dimension $d$ there exists a (d-1)-parameter family
of analytic discs attached to $M$ so that if $f\in C(M)$ holomorphically
extends to those discs, then $f$ is a CR function on $M$.

Despite the large amount of work done on the subject, the conjecture
of Globevnik and Stout has been open so far. In this paper
we prove a version of the conjecture in which the complex lines are
replaced by the complex geodesics of the Kobayashi metric for $D_1$
also known as extremal or stationary discs, whose theory was
developed by Lempert \cite{L81}.
We believe that the extremal discs for a general convex domain $D_1$ are more
appropriate in the problem than the lines because they are intrinsically
defined, invariant under biholomorphisms, and coincide with the lines
for the ball. Hence, if $D_1$ is the ball and $D_2$ is an arbitrary
strictly convex subdomain, then our result proves the conjecture of
Globevnik and Stout for the lines as stated. As in Globevnik's result
\cite{G83} cited above, we only need the
extendibility into the extremal discs tangent to a sufficiently large
open set in $\partial D_2$ (cf. Remark~\ref{r3.1}).

The authors of the results for the concentric balls use the Fourier
analysis and decomposition into spherical harmonics.
This method does not seem to work for general convex domains.
We employ the method of \cite{T04} according to which we add an extra
variable, the fiber coordinate in the projectivized cotangent
bundle. Then using the lifts of the extremal discs we lift the given
function $f$ to a CR function on the boundary of a wedge $W$ whose edge is
the projectivized conormal bundle of $\partial D_1$. Then using
the theory of CR functions we extend it to a bounded holomorphic
function in $W$. Finally since $W$ contains ``large'' discs, we prove that
the lifted function actually does not depend on the extra variable,
which proves the result.

We feel that the method developed here has a wider scope, and we plan
to use it on other occasions.

\section{ Extremal discs}
We will collect here, and develop in some details, the main results of \cite{L81} which are
needed for our discussion.
Let $D$ be a bounded domain of $\C^n$ with $C^k$-boundary; according to \cite{L81} we
assume $k\geq 6$.
We also assume that $D$ is strongly convex in the sense that $D$ has a global defining
function with positive real Hessian.  An analytic disc in $\C^n$ is a holomorphic
mapping $\Delta\to \C^n$, smooth  up to $\partial\Delta$, where $\Delta$ is the standard disc
in $\C$. We denote by $A$ the image set under $\phi$. The disc $A$ is said to be  ``attached''
to $\partial D$ when $\partial A\subset \partial D$.
\bd
\Label{d2.1}
An analytic disc $\phi$ in $D$ is said to be ``stationary'' when it is attached to $\partial D$ and endowed with
a meromorphic lift  $\phi^*(\tau)\in (T^*\C^n)_{\phi(\tau)}\,\,\forall \tau\in\Delta$  with one simple pole at $0$
such that $\phi^*(\tau)\in (T^*_{\partial D}\C^n)_{\phi(\tau)}$ when $|\tau|=1$.
In other words, $(\phi,\phi^*)$ is attached to the conormal bundle $T^*_{\partial D}\C^n$.
\ed
\bd
\Label{d2.2}
An analytic disc $\phi$ in $D$ is said to be ``extremal'' when for any other disc $\psi$ in $D$ with $\psi(0)=\phi(0)$
and $\psi'(0)=\lambda\phi'(0),\,\lambda\in\C$, we have $|\lambda|<1$.
\ed
It is shown in \cite{L81} that extremal and stationary discs coincide. Also, it is shown that they are stable under
reparametrization. In particular,
in Definitions~\ref{d2.1} and \ref{d2.2} we can replace $0$ by any other value of $\tau\in\Delta$ which does not
affect the stationarity or extremality of $\phi$.
It follows  that the extremal discs are the geodesics of the Kobayashi metric in $D$;
in particular they are embeddings of $\bar\Delta$ into $\C^n$.
We recall basic facts about the existence, uniqueness, and smooth dependence
of the extremal discs on parameters (see \cite{L81}, Proposition 11'):
\begin{itemize}
\item[(2.0)]
For any $z\in D$ and $v\in\dot\C^n:=\C^n\setminus\{0\}$, there exists a unique  extremal disc $\phi=\phi_{z,v}$
such that $\phi(0)=z$ and $\phi'(0)=rv$ for $r\in\R^+$. Also, the mapping
\begin{equation*}
D\times \dot \C^n\to C^{2,\frac12}(\bar\Delta),\,\,(z,v)\mapsto \phi_{z,v} \T{ is of class $C^{k-4}$},
\end{equation*}
where $C^{2,\frac12}$ is the space of functions whose derivatives up to order $2$ are $\frac12$-H\"older-continuous.
\end{itemize}
If $\phi^*$ has its pole at $\tau_o$, we multiply it by
$$
\nu(\tau)=\frac{(\tau-\tau_o)(1-\bar\tau_o\tau)}{\tau},\quad\tau\in\Delta,
$$
so that the pole is moved to $0$. Next, we  multiply $\phi^*$ by a real constant $\neq0$ so that $\phi^*(1)$ is the unit outward conormal to $D$ at $\phi(1)$. We will assume that $\phi^*$ is normalized by the two above conditions.
It is essential for our discussion also to clarify the dependence of $\phi^*_{z,v}$  on the parameters $z$ and $v$ which is not explicitly stated in \cite{L81}.
\bp
\Label{p0.1}
The mapping
\begin{equation}
\Label{phi*zv}
D\times \dot\C^n\to C^{2,\frac12}(\bar\Delta),\quad (z,v)\mapsto \phi^*_{z,v},
\end{equation}
is $C^{k-4}$.
\ep
\bpf
Our starting remark is that we can describe $\phi^*_{z,v}$ quite explicitly only over $\partial\Delta$. In fact,  we must have
\begin{equation}
\Label{0.10}
\phi^*_{z,v}(\tau)=g_{z,v}(\tau)\partial\rho(\phi_{z,v}(\tau))\quad\forall \tau\in\partial\Delta,
\end{equation}
where $g_{z,v}$ is real and normalized by the condition $g_{z,v}(1)=1$. On the other hand, when evaluating $\phi^*_{z,v}$ at points $\tau\in\Delta$, we can use the Cauchy integral over $\partial\Delta$. Thus, if we are able to show that $(z,v)\mapsto\phi^*_{z,v}$ with values in $C^{2,\frac12}(\partial\Delta)$ is $C^{k-4}$, the same will be true with values in $C^{2,\frac12}(\bar\Delta)$ since the Cauchy integral preserves  fractional regularity. Now, the second term on the right side of \eqref{0.10} can be handled by means of (2.0) and so what is needed is to describe $g_{z,v}$. We can suppose without loss of generality $\partial_{z_1}\rho\neq0$ for any point of $\phi_{z,v}(\partial \Delta)$. According to \cite{L81} Proposition 9,  we can further normalize our coordinates
in a neighborhood of the disc $\phi_{z,v}(\Delta)$
so that
\begin{equation}
\partial_{z_1}\rho(\phi_{z,v}(\tau))\approx\bar\tau\quad\T{ on }\partial\Delta,
\end{equation}
where ``$\approx$'' means close in $C^1$-norm. In particular, the index of the curve $\left\{ \tau\partial_{z_1}\rho(\phi_{z,v}(\tau))\ \tau\in\partial\Delta \right\}$ around $0$ is $0$ and hence it is well defined the function $\log\left( \tau\partial_{z_1}\rho(\phi_{z,v}(\tau))\right)$ that we will denote by $f$. We are thus reduced to solve the Riemann problem of finding $G=G_{z,v}$ holomorphic such that
\begin{equation}
\Im G_{z,v}(\tau)=\Im\left(\T{log}(\tau\partial_{z_1}\rho(\phi_{z,v}(\tau))\right).
\end{equation}
To this end we have just to set $G=-T_1(\Im (f))+i\Im(f)$ where $T_1$ is the Hilbert transform normalized by the condition $T_1f(1)=1$. Since $T_1$ preserves fractional regularity, then $(z,v)\mapsto G_{z,v}\in C^{{2,\frac12}}(\partial\Delta)$ is also $C^{k-4}$. We finally put
$$
\nu_{z,v}(\tau):=\frac{e^{G(\tau)}}{\tau\partial_{z_1}\rho(\phi_{z,v})}.
$$
We have
\begin{gather}
\begin{split}
\Label{0.11}
\nu_{z,v}&=exp\left(\Re G+i\Im\T{log}\left(\tau\partial_{z_1}\rho(\phi_{z,v})\right)-\T{log}(\tau\partial_{z_1}\rho(\phi_{z,v})\right)
\\
&=exp(\Re G-\Re\T{log}(\tau\partial_{z_1}\rho))\T{ is real },
\end{split}
\\
\nu_{z,v}\tau\partial_{z_1}\rho(\phi_{z,v})\T{ extends holomorphically from } \partial\Delta\T{ to }\bar\Delta.
\end{gather}
Finally, since
$$
\frac g\nu\big|_{\partial\Delta}\in\R,\quad \frac g\nu\T{ extends holomorphically }, \quad \frac g\nu(1)=1,
$$
then $g\equiv \nu$. It follows that $g$ and hence also $\varphi^*$ depends on a $C^{k-4}$ fashion on $z,\ v$.

\epf
There is a statement perfectly analogous to (2.0) above, in which vectors $v$
are replaced by covectors $\zeta$: for any $z\in D$ and $\zeta\in\dot \C^n$ there are unique
the stationary disc and its lift $(\phi,\phi^*)=(\phi_{z,\zeta},\phi^*_{z,\zeta})$ such that
$\phi(0)=z$ and $\phi^*(0)=\zeta$ where $\phi^*(0)$ stands for the residue $\T{Res}\,\phi^*(0)$.

We recall now some basics about the Lempert Riemann mapping.
For any pair of points $(z,w)$ in $D$, let $\phi_{z,w}$ be the (unique) stationary disc through
$z$ and $w$ normalized by the condition $z=\phi_{z,w}(0),\,\,w=\phi_{z,w}(\xi)$ for some $\xi\in(0,1)$;
we define
$$
\Psi_z(w):=\xi\frac{\phi'_{z,w}(0)}{|\phi'_{z,w}(0)|}.
$$
Let $\Bbb B^n$ (resp. $\Bbb D$) denote the unit ball of $\C^n$ (resp. the diagonal of $\C^n\times\C^n$), and set $\dot{\Bbb B}^n=\Bbb B^n\setminus\{0\}$.
Consider the correspondence
\begin{equation}
\Label{Riemann}
(D\times D)\setminus\Bbb D\to D\times\dot{\Bbb B}^n,\quad (z,w)\mapsto (z,\Psi_z(w)).
\end{equation}
 We have
\begin{itemize}
\item
For fixed $z$, $\Psi_z$ is a diffeomeorphism of class $C^{k-4}$ which extends as a diffeomorphism between the boundaries $\partial D$ and $\partial\B^n$.
\item
\eqref{Riemann} is differentiable of class $C^{k-4}$.
\end{itemize}
Write $v=v(z,w)$ for $\Psi_z(w)$. By the above statements, the smoothness of (2.0) and \eqref{phi*zv} are equivalent  to those of
\begin{equation}
\Label{phibis}
(z,w)\mapsto \phi_{z,w},\quad (z,w)\mapsto \phi^*_{z,w}.
\end{equation}
\br
Let $z_\nu$  and $w_\nu$ be sequences converging to $z_o$, and put $v_\nu:=\phi'_{z_\nu,w_\nu}(0)$. If we
 define $v:=\underset\nu{\T{lim}}\,\,\frac{w_\nu-z_\nu}{|w_\nu-z_\nu|}$, then  $v=\underset\nu{\T{lim}}\,\,\frac{v_\nu}{|v_\nu|}$.
 Hence we have convergence (in the $C^{2,\frac12}(\bar\Delta)$ space):
\begin{equation}
\Label{limit}
\phi_{z_\nu,w_\nu}(=\phi_{z_\nu,v_\nu})\to\phi_{z_o,v},\qquad
\phi^*_{z_\nu,w_\nu}(=\phi^*_{z_\nu,v_\nu})\to\phi^*_{z_o,v}.
\end{equation}
\er
For our further needs it is convenient to state the following uniqueness theorem which is largely contained in former
literature.
\bt
\Label{t0.1}
Let us be given two stationary discs $\phi_j\,\,j=1,2$ in a strongly convex domain $D$
and assume that for $\tau_j\in\Delta\,\,j=1,2$ we have
\begin{equation}
\begin{cases}
\phi_1(\tau_1)=\phi_2(\tau_2),
\\
\phi^*_1(\tau_1)=\lambda\phi_2^*(\tau_2)\T{ for some $\lambda\in\C$}.
\end{cases}
\end{equation}
Then, after reparametrization of $\Delta$ we have, for a complex scalar function $\mu=\mu(\tau)$:
$$
\phi_1=\phi_2\quad\phi^*_1=\mu\phi^*_2.
$$
\et
As before, if $\tau_j$ is a pole of $\phi^*_j$, then $\phi^*_j(\tau_j)$ stands for $\T{Res}\,\phi^*_j(\tau_j)$.
\bpf
We assume that the poles of the $\phi^*$'s are placed at $0$.
We compose each $(\phi_j,\phi_j^*)$ with an automorphism of $\Delta$ which brings $\tau_j$ to $0$.
We are therefore reduced to the following:
\begin{equation}
\Label{0.2}
\begin{cases}
\phi_1(0)=\phi_2(0),
\\
\T{Res}\,\phi^*_1(0)=\lambda\T{Res}\,\phi^*_2(0),
\end{cases}
\end{equation}
for a new constant $\lambda$. We put $\lambda=re^{i\theta}$ and replace
$(\phi_2(\tau),\phi_2^*(\tau))$ by  $(\phi_2(e^{-i\theta}\tau),r\phi^*_2(e^{-i\theta}\tau))$.
This transformation reduces \eqref{0.2} to $\lambda=1$. At this point we can prove that $\phi_1=\phi_2$.
We reason by contradiction and suppose $\phi_1\neq\phi_2$. It follows
\begin{equation}
\Label{0.3}
\int_0^{2\pi}\Re{\langle \phi_1^*(\tau)-\phi_2^*(\tau),\phi_2(\tau)-\phi_1(\tau)\rangle}d\theta>0,
\end{equation}
since the integrand is almost everywhere $>0$ on $\partial\Delta$ due to the strong convexity of the domain.
On the other hand $d\theta=-i\frac{d\tau}\tau$; also, $\frac{\phi_2-\phi_1}\tau$ and $\phi^*_1-\phi^*_2$
are holomorphic. Hence the integrand in \eqref{0.3} is a $(1,0)$ form whose coefficient is the real part
of a holomorphic function. Hence the integral \eqref{0.3} is $0$, a contradiction.
\epf
In particular in the situation of Theorem~\ref{t0.1} we have coincidence of the image sets
$\phi_1(\Delta)=\phi_2(\Delta)$.
\br
\Label{r0.1}
Let $\dot T^*\C^n$ be the cotangent bundle to $\C^n$ with the $0$-section removed,
and let $\dot T^*\C^n/\dot\C\simeq\C^n\times\Bbb P^{n-1}_\C$ be the projectivization
of its fibers. We
denote by $(z,[\zeta])$ the variable in $\dot T^*\C^n/\dot \C$.
We can rephrase Theorem~\ref{t0.1} by saying that
if  two discs $(\phi_j,[\phi^*_j])\,\,j=1,2$ have a common point, then, after reparametrization, they need to coincide.
Also, it is useful to point out that, given a stationary disc $\phi(\Delta)$, its lift $[\phi^*(\Delta)]$ is unique.
In fact, the different choices of  $\phi$ obtained by reparametrization, do not affect the class of $\phi^*$ in the
projectivization of the cotangent bundle.
\er

\section{The main result}
Let $D_1$ and $D_2$ be bounded domains of $\C^n$ with $D_2\subset\subset D_1$. We assume that $D_1$
is strongly convex and with $C^k$ boundary for $k\geq 6$ as is the setting of \cite{L81}.
Let $D_2$ be defined by $\rho<0$ for a real function $\rho$ of class $C^2$ with $\partial\rho(z)\neq0$
when $\rho(z)=0$.
\bd
\Label{kobayashiconvex}
The domain $D_2$ is said to be strongly convex with respect to the extremal discs of $D_1$,
if  every such disc $\phi$ tangent to $D_2$ at $z_o=\phi(0)\in\partial D_2$ has tangency of order 2,
that is  for some $c>0$ we have $\rho(\phi(\tau))\ge c|\tau|^2\,\,\forall \tau\in\Delta$, in particular $\bar D_2\cap \phi(\Delta)=\{z_o\}$.
\ed
Here is the main result of our paper.
\bt
\Label{t3.1}
Let  $D_2\subset\subset D_1$
be bounded domains of $\C^n$ with $D_1$ strongly convex and $C^k$ for $k\geq6$, and $D_2$ strongly
convex with respect to the extremal discs of $D_1$ and of class $C^2$.
Let $f$ be a continuous function which extends holomorphically along each extremal disc $\phi(\Delta)$ of $D_1$ which is tangent to $\partial D_2$. Then $f$ extends holomorphically to $D_1$, continuous up to $\partial D_1$.
\et
\br
\Label{}
We do not think that the assumption that $D_2$ is strongly convex with respect to the extremal discs is essential.
We add it for the sake of simplicity and convenience of the proof.
\er
\bpf
We consider the cotangent (respectively tangent) bundle $\dot T^*\C^n/\dot\C$, resp. $\dot T\C^n/\dot\C$, with projectivized fibers $\Bbb P^{n-1}_\C$
and with coordinates $(z,[\zeta])$ and  $(z,[v])$ respectively. The prefix $T^\C$ will be used to denote the complex tangent bundle.
We fix a rule for selecting a ``distinguished'' representative $v$ of $[v]$   and define a mapping
\begin{gather}
\Label{Phi}
(\dot T^\C\partial D_2/\dot\C)\times\Delta\overset\Phi\to(\dot T^*\C^n/\dot \C)|_{D_1\setminus D_2},
\\
\Label{Pphi}
(z,[v],\tau)\overset\Phi\mapsto\left(\phi_{z,v}(\tau),[\phi^*_{z,v}(\tau)]\right),
\end{gather}
where $\phi_{z,v}$ is the unique stationary disc such that $(\phi(0)=z,\,\phi'(0)=rv)$
for some $r\in\R^+$ and $\phi^*_{z,v}$ is its ``lift''
according to \S 2, (2.0).
Note that by multiplying $\phi^*_{z,v}$ by $\nu(\tau)$, real on $\partial\Delta$, we can move the pole to $\tau=0$.

We denote by $\mathcal S$ the image-set of $\Phi$.
 We show  that $\Phi$ is  an injective smooth local parametrization of $\mathcal S$.  First, it is injective:
 in fact, if $(z_1,[v_1],\tau_1)$ and $(z_2,[v_2],\tau_2)$ go to the same image, then by Theorem~\ref{t0.1}
 in \S 2, $\phi_{z_1,v_1}$ and $\phi_{z_2,v_2}$ coicide up to reparametrization. On the other hand, by the strong
 convexity of $D_2$ with respect to the stationary discs of $D_1$, we must have $z_2=z_1$. Then we also have $v_1=v_2$
 by our rule of taking representatives and therefore  the discs coincide (without need of reparametrization).
 Finally  $\tau_2=\tau_1$ because they are injective (cf. \S 2).
As for the smoothness, we make a  choice of our representative $v$ smoothly depending on $z$, and point our
attention to (2.0) and \eqref{phi*zv} of \S 2. If we then take evaluation of the discs and their lifts at
$\tau\in\Delta$ we get the $C^{k-4}$-smoothness of \eqref{Pphi}. In the lines of what was remarked after
Proposition~\ref{p0.1}, for any $(z,[\zeta])\in(\dot T^*\C^n/\dot C))|_{D_1\setminus D_2}$ there is a
unique $(\phi,\phi^*)$, up to reparametrization, such that $\phi(\tau)=z,\,\,[\phi^*(\tau)]=[\zeta]$ for some
$\tau\in\Delta$. On the other hand, the class of stationary discs which are tangent to $\partial D_2$ divides
the set of all stationary discs into two sets, the ones which are transversal to (resp. disjoint from) $D_2$.
Accordingly, $\mathcal S$ divides $(\dot T^*\C^n/\dot\C)|_{D_1\setminus\bar D_2}$ into two sets.
We denote by $\mathcal W$ the first set and refer to it as to the ``finite" side of $\mathcal S$ the complement
being called a neighborhood of the ``plane at infinity". The set
$\mathcal W$ is a wedge type domain with boundary $\mathcal S$ and edge
$\mathcal E:=\dot T^*_{\partial D_1}\C^n/\dot\C$.

We now describe the fibers $\mathcal S_{z_o}=\pi^{-1}(z_o)\cap\mathcal S$ where $\pi:\dot T^*\C^n/\dot\C\to\C^n$
is the projection $\pi(z,[\zeta])=z$.
Our plan is to use $\Psi_{z_o}$, interchange  $D_1$ with $\B^n$ and $z_o$ with $0$, analyze the situation in this
new setting, and then bring back the conclusions to the former by $\Psi^{-1}_{z_o}$.
Recall that $\Psi_{z_o}$ interchanges the stationary discs through $z_o$ with the complex lines (the stationary
discs of the ball) through $0$. We first describe the set
\begin{equation}
\Label{gamma}
\gamma_0=\{z\in\partial (\Psi_{z_o}D_2):\,\,\T{ for some }v\in T^\C_z\partial(\Psi_{z_o} D_2),\,\,
\phi_{z,v}\T{ passes through }0\}.
\end{equation}
If $\rho=0$ is an equation for $\partial(\Psi_{z_o} D_2)$, $\gamma_0$ is  defined by
$$
\rho(z)=0\quad\partial\rho(z)\cdot z=0.
$$
This is a system of three real equations that we denote by $r=0$.
We  normalize our coordinates so that
$$
\quad\partial\rho(z)=(1,0,\dots),\quad z=(0,c,0,\dots).
$$
We  then have for the partial Jacobian $J_{z_1,\bar z_1,z_2,\bar z_2}r(z)$:
\begin{equation}
\Label{det1}
\left[
\begin{matrix}
1&1&0&0
\\
&&&
\\
* & * &c\partial^2_{z_2,z_2}\rho&c\partial^2_{\bar z_2,z_2}\rho
\\
&&&
\\
* & * & c\partial^2_{z_2,\bar z_2}\rho& c\partial^2_{\bar z_2,\bar
 z_2}\rho
\end{matrix}
\right],
\end{equation}
where the asterisks denote unimportant matrix coefficients.
Let $A$ be the $3\times 3$ minor obtained by discarding the first column. We have
\begin{equation}
\Label{det2}
\begin{split}
\T{det}\,A &=c^2\T{det}\left[
\begin{matrix}
\partial^2_{z_2,z_2}\rho&\partial^2_{z_2,\bar z_2}\rho
\\
\partial^2_{\bar z_2,z_2}\rho & \partial^2_{\bar z_2,\bar z_2}\rho
\end{matrix}
\right]
\\
&=-c^2\T{det}\,\left(Hess (\rho)|_{\C z_2}\right)\,\,<\,\,0,
\end{split}
\end{equation}
where the real Hessian of $\rho$ at $z$ along the $z_2$-plane is positive
because $\Psi_{z_o}D_2$ is strongly convex with respect to $\Psi_{z_o}\left(\phi_{z_o,z}(\Delta)\right)$.
In conclusion, $\T{rank}\,J(r) =3$ and hence $\gamma_0$ is a regular real manifold of dimension $2n-3$  compact
and without boundary. We now use the fact that $\Psi_{z_o}$ is a diffeomeorphism and conclude that
$\gamma_{z_o}:=\Psi^{-1}_{z_o}(\gamma_0)$ is also a regular manifold of dimension $2n-3$ in $\partial D_2$,
which enjoys the same properties as $\gamma_0$. It represents the set of points where the geodesics of $D_1$
through $z_o$ are tangent to $\partial D_2$.
Let $\tilde \gamma_{z_o}$ be the section $(z,[v(z)])$ of $(\dot T^\C\partial D_2/\dot\C)\big|_{\gamma_{z_o}}$
where $[v(z)]$ is the direction tangent at $z$ to the stationary disc connecting $z_o$ and $z$.
We can parametrize  the fiber $\mathcal S_{z_o}$ over $\tilde \gamma_{z_o}\times \Delta$ by the same
parametrization $\Phi$ as in \eqref{Phi}. This being bijective, we conclude that $\mathcal S_{z_o}$ is a finite
family of regular closed manifolds of dimension $2n-3$ without boundary, which do not intersect.
We move now $z$ from the fixed $z_o$ and describe the behavior of the fibers $\mathcal S_z$; they depend
in a ($C^{k-4}$) fashion on $z$ since the mapping in \eqref{Riemann} is also $C^{k-4}$.
As for their behavior at $z_o\in\partial D_2$, we consider the set $\Pi_{z_o}$ defined by the diagram
\begin{equation*}
\begin{matrix}
\dot T^\C_{z_o}\C^n/\dot \C&\underset\sim\to&\dot T^*_{z_o}\C^n/\dot\C&\simeq&\Bbb P^{n-1}_\C
\\
&&&&
\\
\cup&&\cup&&
\\
&&&&
\\
\dot T^\C_{z_o}\partial D_2/\dot \C&\underset\sim\to&\Pi_{z_o}&&
\end{matrix}
\end{equation*}
where the two horizontal arrows are given by  the smooth injective mapping $v\mapsto [\phi^*_{z_o,v}(0)]$.
Thus $\Pi_{z_o}:=\{[\phi^*(0)]:\,\, \phi\T{ is tangent to $\partial D_2$ at $z_o$}\}$  is a 2-codimensional
real submanifold of $\Bbb P^{n-1}_\C$ which reduces to a single point when $n=2$.
\bl
The sets $\mathcal S_{z_\nu}$ shrink to $\Pi_{z_o}$ as
$z_{\nu}\to z_o\in\partial D_2$; in particular, $\mathcal S_{z_\nu}$ consists of just one component when $z_\nu$ is close to $z_o$.
\el
\bpf
By the strong convexity of $\partial D_2$, the manifolds $\gamma_{z_\nu}$ shrink to $\{z_o\}$ as $z_\nu\to z_o$.
If we pick up any sequence $w_\nu\in\gamma_{z_\nu}$, we have
$$
\frac{w_\nu-z_o}{|w_\nu-z_o|}\to v\in T^\C_{z_o}\partial D_2.
$$
Let $\phi_{z_\nu,w_\nu}$ (resp. $\phi_{z_o,v}$) be the geodesic through $z_\nu$ and $w_\nu$ (resp. through $z_o$
with tangent $v$),  normalized by the condition
$z_\nu=\phi_{z_\nu,w_\nu}(0),\,\,w_\nu=\phi_{z_\nu,w_\nu}(\xi)$ for $\xi\in(0,1)$,
(resp. $z_\nu=\phi_{z_o,v}(0),\,\,rv=\phi'_{z_o,v}(0)$ for $r\in\R^+$).
Then
$$
\phi_{z_\nu,w_\nu}\to \phi_{z_o,v},\quad \phi^*_{z_\nu,w_\nu}\to\phi^*_{z_o,v},
$$
with convergence in the $C^{{2,\frac12}}(\bar\Delta)$ norm. In particular, since
$\mathcal S_{z_\nu}=\underset{w_\nu\in\gamma_{z_\nu}}\bigcup[\phi^*_{z_\nu,w_\nu}(0)]$, then
$$
\mathcal S _{z_\nu}\to\underset{v\in\dot T^\C_{z_o}\partial D_2}\bigcup[\phi^*_{z_o,v}(0)].
$$

\epf
It follows that for the fibers $\mathcal W_{z_\nu}$, which are open domains of $\Bbb P^{n-1}_\C$ with boundary
$\mathcal S_{z_\nu}$, we have merely by definition:
$$
\mathcal W_{z_\nu}\to \Bbb P^{n-1}_\C\setminus\Pi_{z_o}\T{ as }z_\nu\to z_o\in\partial D_2.
$$
If, instead, we move $z_\nu$ towards $z\in\partial D_1$,
then each $\mathcal S_{z_\nu}$ as well as their ``finite'' sides $\mathcal W_{z_\nu}$, shrink to the single
point $\dot T^*_{\partial D_1}\C^n/\dot\C\big|_z$.

Now we move $z_o$ all over  $ D_1\setminus\bar D_2$. If we  take a closer look to \eqref{det1}, \eqref{det2} we see that the set $\Psi_{z_o}(D_2)$ as well as its equation $\rho_{z_o}=0$, moves smoothly with respect to $z_o$ by the regularity properties of $\Psi$. It follows that the set defined by
$\{(z_o,z):\,\,z_o\in D_1\setminus\bar D_2,\,\,z\in\gamma_{z_o}\}$ is a $(4n-3)$-dimensional manifold. In particular,  the set $\gamma_{z_o}$ is a $(2n-3)$-dimensional  manifold and it cannot turn from one to several components without passing through a singular point $z_o$.
It follows that the set $\mathcal S_{z_o}$ also  consists of one component.

By the preceding discussion and Sard's theorem, we can also say  that $\mathcal S$ is a smooth regular manifold
except possibly a closed subset of measure zero. Along with its natural foliation by the discs
$(\phi_{z,v},[\phi^*_{z,v}])$,
we need to endow $\mathcal S$ with another foliation, locally on a neighborhood of each of its points,
by CR manifolds $\mathcal M$ of dimension $2n$ and CR dimension $1$ each one being still a union of discs.
For this, we fix $z\in \bar D_1\setminus \bar D_2$, consider the submanifold $\gamma_z$ of $\partial D_2$ with
dimension $2n-3$ of  points of tangency for the stationary discs through $z$, and denote by $w$ the point which
moves in $\gamma_z$.
As above, we denote by $\phi_{w,z}$ the stationary disc through $w$ and $z$, normalized by
$\phi_{w,z}(0)=w,\,\,\phi_{w,z}(\xi)=z$ for $\xi\in (0,1)$; we also write $\xi=\xi(w,z)$ and define
$v(w,z):=\phi_{w,z}(0)$.
We set $\Gamma_z:=\{(w,[v(w,z)],\xi(w,z))\,:\, w\in\gamma_z\}$; then $\T{dim}\,\Gamma_z=\T{dim}\,\,\gamma_z=2n-3$.
Since $\Phi_1$ sends all points of $\Gamma_z$ to the fixed $z$, then we have an inclusion
$T\Gamma_z\subset \T{Ker}\,\Phi'_1\big|_{\Gamma_z}$. But since the dimensions  are the same, then
$T\Gamma_z=\T{Ker}\,\,\Phi'_1\big|_{\Gamma_z}$.
In particular,
if $p$ is the projection $(w,[v],\tau)\mapsto w$, then
\begin{equation}
\Label{Ker}
p'\left(\T{Ker}\,\,\Phi'_1\big|_{\Gamma_z}\right)\subset T\gamma_z.
\end{equation}
We define $\mathcal M$ locally at a point $(z,[\zeta])\in\mathcal S\cup\mathcal E$; if
$(z,[\zeta])\in\mathcal E$, $\mathcal M$ will be in fact a manifold with boundary $\mathcal E$.
Let $(w,[v],\tau)$ be the value
 of the parameter in $(T^\C\partial D_2)\times\Delta$ which corresponds to $(z,[\zeta])$ via $\Phi$.
Choose a germ   of submanifold  $\delta_z\subset\partial D_2$ transversal to $\gamma_z$ at $w$ with
complementary dimension $2$.
By \eqref{Ker}, we have
\begin{equation}
\T{Ker}\,\,\left(\Phi'_1(w,[v],\tau)\big|_{T_w\delta_z\times\Bbb P^{n-1}_\C\times T_\tau\Delta}\right)=\{0\}.
\end{equation}
Thus $\Phi_1$ induces a diffeomeorphism between a neighborhood $\Sigma=\Sigma_1\times\Sigma_2$ of $(w,[v],\tau)$
in $(T^\C\partial D_2/\dot\C)|_{\delta_z}\times\Delta$ and a neighborhood of $z$ in $\bar D_1$.
We define $\mathcal M=\Phi(\Sigma)$ that is
\begin{equation}
\Label{2.6}
\mathcal M=\underset {(\phi,\phi^*)}\cup(\phi,[\phi^*])(\Sigma_2),
\end{equation}
for $(\phi(0),[\phi'(0)])\in\Sigma_1$. $\Phi$ is a diffeomorphic parametrization of $\mathcal M$ over $\Sigma$ and
hence $\mathcal M$ is a smooth manifold, in fact a graph over a neighborhood of $z$ in $\bar D_1$.
This was not necessarily the case of $\mathcal S$ since $\Phi$ is a smooth and bijective parametrization
of $\mathcal S$ but it might occur that $\Phi'$ is degenerate at some point.
We define a function $F$ on $\mathcal S$ by collecting all extensions $f_{\phi(\bar\Delta)}$ of the given $f$
from $\phi(\partial\Delta)$ to $\phi(\bar\Delta)$.
For $(z,[\zeta])\in\mathcal S$ we put
\begin{equation}
\Label{main}
F(z,[\zeta])=f_{\phi(\bar\Delta)}(z)\T{ if $(\phi(\tau),[\phi^*(\tau)])=(z,[\zeta])$ for some $\tau$}.
\end{equation}
According to Theorem~\ref{t0.1}, $F$ is well defined.
We have the following
\bp
\Label{p2.1}
At every point of $\mathcal S\setminus \mathcal E$, the function $F$ holomorphically extends to a one-sided
neighborhood on the $\mathcal W$-side of $\mathcal S$.
\ep
\bpf
The ingredients of the proof are the foliation of $\mathcal S$ by manifolds with boundary $\mathcal M$,
which are themselves union of discs, and the additional transversal foliation of $\mathcal W$ by the fibers
$\mathcal W_z$. The starting remark is that $F$ is holomorphic along each disc and therefore it is CR on each
$\mathcal M$ since $\T{dim}_{CR}\mathcal M=1$.

\noindent
{\bf (a)}
We approximate  $F|_{\mathcal E}$ by a sequence of entire functions $F_\nu$ (cf. e.g. \cite{B91}).
To this end it is important to notice, as it was first pointed out by Webster, that since $\partial D_1$ is
strongly convex, then $\mathcal E$ is totally real. Then, in an identification $\mathcal E\simeq\R^{2n-1}$,
these are defined by
\begin{equation}
\Label{approximation}
\hat F_\nu(\xi)=\left(\frac\nu\pi\right)^{\frac{2n-1}2}\int_{\R^{2n-1}}F(\eta)e^{-\nu(\eta-\xi)^2}dV,
\end{equation}
($dV$ being the element of volume in $\R^{2n-1}$).
It is well known that $\hat F_\nu\to F$ uniformly on compact subsets of $\mathcal E$. Also, $F$ being CR
on each $\mathcal M$, it is possible to deform the integration chain from $\mathcal E$ to another chain entering
inside $\mathcal M$ and reaching any point of $\mathcal M$ in a neighborhood of $\mathcal E$. In other terms,
the function $F$ is approximated, over each $\mathcal M$  near $\mathcal E$, by the same sequence \eqref{approximation}
of entire function. Since the $\mathcal M$'s give a foliation of $\mathcal S$,  it follows that the uniform
approximation of $F$ by the $F_\nu$'s holds on the whole $\mathcal S$ in a neighborhood of $\mathcal E$.

\noindent
{\bf (b)} By using now the foliation of $\mathcal W$ by the fibers $\mathcal W_z$, we can bring the approximation
by entire functions from $\mathcal S$ to $\mathcal W$ in a neighborhood of $\mathcal E$: in fact, by maximum principle,
the sequence $\hat F_\nu$ which is Cauchy over $\mathcal S_z$ will be Cauchy on the whole ${\mathcal W}_z$.

\noindent
{\bf (c)}
We use now the theory of propagation of wedge extendibility along discs for each CR function
$F|_{\mathcal M}$ by \cite{T94} which develops \cite{HT81}. We put a suffix $s$ in the notation of the disc
$\Delta_s$ to specify its radius, and define
\begin{equation}
\Label{maximal}
I=\{r\in(0,1):\,\,F\T{ extends to the side $\mathcal W$ of
$\mathcal S$ in }\underset {(\phi,\phi^*)}\cup(\phi,[\phi^*])(\Delta_1\setminus\bar\Delta_{1-r})\},
\end{equation}
for all stationary discs $\phi$ tangent to $\partial D_2$ at $\tau=0$. (The last requirement is just a choice
of the parametrization.) We have $I\neq\emptyset$ due to (b) above.
We show now that we have indeed $I=(0,1)$ from which the proposition follows. We reason by contradiction,
suppose $I\neq(0,1)$, and denote by $r_o$ the supremum of $I$; thus $r_o<1$. By propagation of wedge extendibility
of $F|_{\mathcal M}$ for each $\mathcal M$, and since the wedge evolves continuously with the base point, then,
on account also of a compactness argument, $F$ would extend to the side $\mathcal W$ for a value of $r$ bigger
than $r_o$, a contradiction.
\epf
\noindent
{\it End of proof of Theorem 3.2.}

\begin{itemize}
\item
First, recall that for $z$ moving from $\partial D_1$ to $z_o\in\partial D_2$, the fibers ${\mathcal W}_z$
grow from a single point to ${\mathcal W}_{z_o}=\Bbb P^{n-1}_\C\setminus\Pi_{z_o}$.
Also, recall that by approximation, $F$ extends holomorphically from $\mathcal S_z$ to ${\mathcal W}_z$
when $z$ is close to $\partial D_1$,
and, by propagation, to a neighborhood of $\mathcal S_z$ in $\mathcal W$ when $z$ is no longer
close to $\partial D_1$.
Then  $F$ extends to the whole set ${\mathcal W}$ by the Hartogs continuity principle.
For $n>2$ the same conclusion also follows by the Hartogs extension theorem.
\item
The boundary values of $F$ on
$\pi^{-1}(\partial D_2)\cap\bar{\mathcal W}\subset\partial{\mathcal W}$ are constant on the
fibers ${\mathcal W}_{z_o}$, $z_o\in\partial D_2$.
Indeed, $F|_{{\mathcal W}_{z_o}}$ holomorphically extends to the whole projective space $\Bbb P^{n-1}_\C$
because the set $\Pi_{z_o}$ of codimension 2 is removable, hence it is constant.
Now since $F$ is constant on the fibers of $\bar{\mathcal W}$ on an open set of the
boundary of ${\mathcal W}$, then $F$ is constant on the fibers of ${\mathcal W}$ everywhere in ${\mathcal W}$.
Then $\tilde f(z):= F(z,[\zeta])$, $(z,[\zeta])\in {\mathcal W}$ is a well defined holomorphic extension of
the original function $f$ to $D_1\setminus\bar D_2$. Then $f$ further extends to $D_2$ by the Hartogs theorem.
The proof is now complete.
\end{itemize}

\epf
\br\Label{r3.1}
Take a line segment $I$ connecting a pair of points $z_1$ and $z_2$ of $\partial D_1$ and $\partial D_2$ resp.,
fix a neighborhood $U\supset I$ in $\C^n$, denote by $\mathcal I$ the family of discs tangent to $\partial D_2$
and passing through $U\cap D_1$, and set $V:=\underset{\phi\in\mathcal I}\cup\phi(\partial\Delta)$.
Assume that $f$ is defined and continuous in $V$ and extends holomorphically to the discs which belong
to $\mathcal I$; then $f$ extends holomorphically to a one-sided neighborhood of $V$ in $D_1$.
In fact, by moving $z$ from $z_1$ to $z_2$ along $I$ we will have the same conclusions for the fibers
$\mathcal S_z$ and $\mathcal W_z$ as in the proof of Theorem~\ref{t3.1}. In particular we will conclude that $F$
is independent of $[\zeta]$ in a neighborhood of $z_1$. But then $F$ is independent of $[\zeta]$ wherever
it is defined, in particular in a one-sided neighborhood of $V$.
\er


\begin{thebibliography}{BER99}

\bibitem[1]{AS}{\bf M. L. Agranovsky, A. M. Semenov}--- Boundary analogues
of the Hartogs theorem. {\em Sibirsk. Mat. Zh.} {\bf  32}  (1991),  no. 1,
168--170, 222;  translation in  Siberian Math. J.  32  (1991),  no. 1,
137--139

\bibitem[2]{BCPZ}{\bf  C. Berenstein, D. C. Chang, D. Pascuas,
L. Zalcman}---Variations on
the theorem of Morera. {\em The Madison Symposium on Complex Analysis}
(Madison, WI, 1991),  63--78, Contemp. Math., {\bf 137}, Amer. Math. Soc.,
Providence, RI, (1992)

\bibitem[3]{B91} {\bf  A. Boggess} --- CR manifolds and the
tangential
Cauchy-Riemann complex. {\em Studies in Adv. Math. CRC Press},
(1991)

\bibitem[4]{GS}{\bf J. Globevnik, E. L. Stout}--- Boundary Morera theorems
for holomorphic functions of several complex variables. {\em Duke Math. J.}
{\bf 64}  (1991),  no. 3, 571--615

\bibitem[5]{G83} {\bf J. Globevnik}--- On holomorphic extensions from
spheres in $C\sp{2}$. {\em Proc. Roy. Soc. Edinburgh} Sect. A  {\bf 94}
(1983), no. 1-2, 113--120

\bibitem[6]{G87}{\bf J. Globevnik}--- A family of lines for testing
holomorphy in the ball of $C\sp 2$. {\em Indiana Univ. Math. J.} {\bf  36}
(1987),  no. 3, 639--644

\bibitem[7]{HT81}{\bf J. Hanges, F. Treves}---Propagation of holomorphic
extendability of CR functions.  {\em Math. Ann.} {\bf 263}  (1983),  no.
2, 157--177

\bibitem[8]{L81}{\bf L. Lempert}---La m\'etrique de Kobayashi et la
representation des domaines sur la boule.
 {\em Bull. Soc. Math. de France} {\bf 109} (1981),
427--474

\bibitem[9]{R}{\bf W. Rudin}--- Function theory in the unit ball of
$C\sp{n}$. {\em Grundlehren der Mathematischen Wissenschaften}
[Fundamental Principles of Mathematical Science], {\bf 241}.
Springer-Verlag, New York-Berlin, (1980)

\bibitem[10]{S77} {\bf E. L. Stout}--- The boundary values of holomorphic
functions of several complex variables. {\em Duke Math. J.} {\bf 44}
(1977), no. 1, 105--108


\bibitem[11]{T94} {\bf A. Tumanov}---Connections and propagation of
analyticity for CR functions. {\em Duke Math. J.} {\bf 73}  (1994),  no.
1, 1--24


\bibitem[12]{T97}{\bf A. Tumanov}---Thin discs and a Morera theorem for CR
functions.  {\em Math. Z.}  {\bf 226}  (1997),  no. 2, 327--334


\bibitem[13]{T04}{\bf A. Tumanov}---A Morera type theorem in the strip.
{\em Math. Res. Lett.}  {\bf 11}  (2004),  no. 1, 23--29


\end{thebibliography}
\end{document}